\documentclass[12pt]{article}
 
\usepackage{amssymb,amsmath,amsfonts,amssymb, verbatim}
\newcommand{\be}{\begin{equation}}
\newcommand{\ee}{\end{equation}}
\newcommand{\la}{\label}
\newcommand{\ba}{\begin{array}{c}}
\newcommand{\ea}{\end{array}}

\newtheorem{thm}{Theorem}

\newcommand{\del}{\delta}
\newcommand{\eps}{\varepsilon}
\newcommand{\calE}{\mathcal{E}}
\newcommand{\one}{{\mathbf 1}}

\title{On the zero temperature limit of interacting corpora}
\author{Peter Constantin\\Department of Mathematics\\The University of Chicago\\and\\
Andrej Zlato\v s\\Department of Mathematics\\The University of Chicago}

\begin{document}
\maketitle 

\begin{abstract}{We characterize the zero-temperature limits of minimal free energy states for interacting corpora -- that is, for objects with finitely many degrees of freedom, such as articulated rods. These limits are measures supported on zero-level-sets of the interaction potential. We describe a selection mechanism for the limits that is mediated by evanescent entropic contributions.}\end{abstract}

\noindent{\bf{MSC 2000 classification numbers}} 82D30, 58D99 


\section{Introduction}
One of the simplest examples of a corpus is a rod-like molecule. A rod-like molecule is modeled by a director $p$, a vector of length $1$ in space. When a melt (ensemble) of such molecules is in statistical equilibrium, the state of the system is described by a probability distribution. The configuration space of all molecules $p$ is $M= {\mathbb S}^2$, the unit sphere in
${\mathbb R}^3$. The state of the system is a probability measure $\nu = fd\mu$ with $f\ge 0$ and $\int_M fd\mu =1$, where $f$ is the density with respect to the normalized volume element $d\mu$ on the sphere. The ground state of the system is obtained by minimizing the {\it free energy} ${\mathcal E}_b[f]$. The free energy is composed of two terms, an entropic term, representing the tendency of the system to occupy as many degrees of freedom as possible, and an excluded volume interaction, expressing the repulsive interaction between the rod-like particles \cite{ons}:
$$
{\mathcal E}_b[f] = \int_M f\log f d\mu + \frac{b}{2}\int_M U[f]fd\mu.
$$
The potential $U[f]$ is computed from the state of the system, and it is given by the integral
$$
U[f](p) = \int_M k(p,q) f(q)d\mu(q)
$$
where the kernel $k$, in mean-field fashion accounts for the microscopic interaction between molecules. The single  parameter $b\ge 0$ embodies both the intensity of the interaction and the inverse of temperature. The free energy is not convex. As the temperature is lowered, or the intensity of interaction is increased, phase transitions do occur. The isotropic state $f_0\equiv 1$ is the high temperature limit ($b\to 0$). At positive $b$ the minimizers solve a nonlinear nonlocal  Euler-Lagrange equation which one of us termed ``Onsager's equation'' \cite{c-nfp}:
$$
f = \left (Z_b[f]\right)^{-1}\exp{\left (-bU[f]\right)}.
$$ 
(This is different from the ``Onsager equation'' for the dielectric polarization.)
The high temperature asymptotic behavior ($b\to\infty$) of solutions with the Maier-Saupe potential is either isotropic $f=1$, oblate $f = \delta_C$ (where $C$ is a great circle), or prolate $f= \delta_P$ (where $P$ is a point) \cite{ckt,fatkullin,forest,lzz,liuzz}.  However, the prolate state is the only limit of absolute minimizers. 

In this paper we investigate the selection mechanism for the zero temperature limit for interacting corpora. The term ``corpus'' refers to an object with finitely many degrees of freedom,  such as an assembly of articulated rods, or the sticks-and-balls model of a molecule. Corpora may be made of parts, freely articulated or not, themselves simpler corpora. Many problems in nonlinear science have at their core the study of statistics of melts of corpora, under the assumptions of disorder and frustration.  The collection of all corpora in a particular problem 
is the configuration space $M$. This can be rather complicated, even if the corpora are just $n$-gons with unequal sides. Because we are ultimately interested in the situation when each corpus is a rather complex network or graph, it is useful to phrase the equilibrium and kinetic problems broadly, in quite general configuration spaces. In a previous paper \cite{c1} we established the existence of solutions of Onsager's equation in metric spaces and discussed the limit $b\to\infty$. In this paper we describe this limit more precisely and establish a selection criterion. Under quite general conditions, minimizers of the free energy solve Onsager's equation at fixed $b$. As $b\to\infty$ the limit points of minimizers are measures concentrated on sets $A\subseteq M$ with $k(p,q) = 0$ for all $p,q\in A$.  The corpora in the support of any limit measure, are the {\it ur-corpora}, the prototype or ancestor corpora. A selection principle among these is in effect. The selection principle is a consequence of the entropic contribution to the free energy, although this contribution becomes (relatively) vanishingly small compared to the particle interaction. The principle essentially says that if a $k$-non-increasing $\mu$-increasing transformation exists from a neighborhood of $p$ to a neighborhood of some $q\in M$, then $p$ cannot be an ur-corpus. 

The Onsager equation can be phrased if the following data are given: 
a compact metric space $M$ with distance $d$, representing the configuration space, a Borel probability measure $\mu$ on $M$ representing the isotropic state, and a
symmetric nonnegative, bi-Lipschitz interaction kernel $k:M\times M\to {\mathbb R}_+$ representing the corpora interaction. 
If $p= (p_1,\dots p_N)\in M= M_1\times\dots \times  M_N$, if $d\mu(p) = d\mu_1(p_1)\times \dots\times d\mu_N(p_N)$ and if $k(p,q) = \sum_{1}^N k_{j}(p_j, q_j)$ then the Onsager equation factorizes, and the solution, for any $b$, is a product
$fd\mu = \Pi_{1}^N f_jd\mu_j$ of appropriate solutions of the Onsager equation on $M_j$. This is the situation in which the corpora comprise of ``freely articulated'' parts: the states are product measures and are made of independent, non-interacting parts. The zero temperature limits are then products of zero-temperature limits of the parts. We will give
simple examples below of corpora made from several interacting parts, compute the limits and explain them using the selection principle.

\section{Framework and Main Results} 

We consider a compact  separable metric space $M$ with
metric $d$, a Borel probability measure $\mu$ on $M$ and a uniformly bi-Lipschitz function
\be
k:M\times M\to {\mathbb R}_+
\la{kernel}
\ee
obeying
\be
\left\{
\ba
k(p,q) = k(q,p) \quad {\mbox{ for all}}\; p, q \in M,\\
\exists \, L>0 : 
\left | k(p,s)- k(q,s)\right |\le Ld(p,q)\quad {\mbox{for all}}\; p,q,s\in M.
\ea
\right.
\la{kprop}
\ee
For Borel measurable functions $f\ge 0$ with integral equal to one,
\be
\int_Mfd\mu = 1,
\la {norma}
\ee
we consider 
the potential
\be 
U[f](p) = 
\int_M k(p,q)f(q)d\mu(q).
\la{uf}
\ee
We take a  parameter $b>0$ and define the free energy
\be
{\calE}_b[f] = \int_M \left\{\log f + \frac{b}{2}U[f]\right\}fd\mu.
\la{fe}
\ee

\noindent{\bf{Remark.}} In view of the compactness of M, the theory does not change if we add a constant to the kernel $k$; that is why there is no loss of generality in assuming $k\ge 0$. Likewise, because $f(p)f(q)d\mu(p)d\mu(q)$ is symmetric in $p,q$, there is no loss of generality in the assumption $k(q,p) = k(p,q)$. Indeed, for an arbitrary kernel $K$, we can write $K(p,q) = k(p,q) + a(p,q)$ with $k(p,q) = k(q,p)$ and $a(p,q) = -a(q,p)$, and only $k(p,q)$ contributes to ${\calE}_b[f]$.  

\smallskip

We start with the proof of the existence of minimizers of the free energy.

\begin{thm} \la{T.1}
Let $b>0$ be fixed. Let $M$ be compact, metric. Let $k$ satisfy
(\ref{kprop}). Let $\mu$ be a Borel probability on $M$. Then there exists a 
function $g\ge 0$ with $\int_M gd\mu =1$ which achieves the minimum
\be
{\calE}_b[g] = \min_{f\ge 0,\; \int_Mfd\mu =1}{\calE}_b[f].
\la{mini}
\ee
Moreover, $g$ is strictly positive, Lipschitz continuous, and solves the Onsager equation
\be
g = \left(Z_b[g]\right)^{-1}e^{-bU[g]}
\la{geq}
\ee
with
\be
Z_b[g] = \int_Me^{-bU[g](p)}d\mu(p).
\la{zeq}
\ee
\end{thm}

\noindent{\bf Proof.} The proof follows along the lines of the proof of a similar result in {\cite{c1}}. We sketch it here for completeness. We note that if
$f\ge 0$ and $\int_Mfd\mu =1$ then the function $U[f]$ obeys
\be
\left\{
\ba
0\le U[f](p) \le \|k\|_\infty\\
\left | U[f](p) -U[f](q)\right | \le Ld(p,q)
\ea
\right. \la{uprop}
\ee
with $L>0$ given in \eqref{kprop}. We also note that from Jensen's inequality, $\mu(M)=1$, and the normalization $\int_M f d\mu =1$ it follows that 
$$
\int_M f\log fd\mu \ge 0,
$$
and consequently the free energy is bounded below:
$$
{\calE}_b[f] \ge 0.
$$
We take then 
a minimizing sequence $f_j$, 
$$
a = \inf_{f\ge 0, \int_M fd\mu =1}{\calE}_b[f] = \lim_{j\to\infty}{\calE}_b[f_j].
$$
Without loss of generality, by passing to a subsequence and relabeling, we may assume that the measures $f_jd\mu$ converge weakly to a measure $d\nu$. Also, using (\ref{uprop}) and the Arzela-Ascoli theorem, we may pass to a subsequence which we relabel again $f_{j}$, so that $U_{j}= U[f_j]$ converge uniformly to a non-negative Lipschitz continuous function $\bar U$. Then it follows that 
$$
\bar U(p) = \int_M k(p,q)d\nu(q)
$$
holds and 
$$
\lim_{j\to\infty}\int_M U_j f_j d\mu = \int_M \bar Ud\nu.
$$
Because ${\calE}_b[f_j]$ is a convergent sequence, the sequence $\int_Mf_j\log f_j d\mu$
also converges and consequently the integrals $\int_M f_j\log_+ f_j d\mu$ are uniformly bounded. 

It then follows that $d\nu$ is absolutely continuous, that is, $d\nu = gd\mu$ 
with $g\ge 0$ and $g\in L^1(d\mu)$. Indeed, this is because the sequence $f_jd\mu$ is uniformly absolutely continuous. The latter is proved using the convexity of the function $y\log y$ and Jensen's inequality
$$
\mu(A)^{-1}\int_A f_j\log f_j d\mu \ge m\log m,
$$
where $m = m(A,j) = \mu(A)^{-1}\int_A f_jd\mu$. Thus for all $A,j$,
$$
m\log m\le \frac C{\mu(A)}
$$
with a fixed  $C\ge 1$. 
Let us choose $R=R(A)\ge 1$ so that $R\log R = C/\mu(A)$. Then $m\le R$ and so $\int_A f_j d\mu \le \mu(A) R =  C/\log R$. This means 
$$
\int_A f_jd\mu \le \delta(\mu(A))
$$
with $\lim_{x\to 0}\delta(x) = 0$ and $\delta(\cdot)$ independent of $A$ and $j$. It follows that $\nu$ is absolutely continuous with respect to $d\mu$.

The weak convergence tested on the function $1$ gives $\int_M gd\mu =1$. In general, weak convergence of measures is not enough to show lower semicontinuity of nonlinear integrals or almost everywhere convergence. We claim however that, in fact, the convergence $f_n\to g$ takes place strongly in $L^1(d\mu)$:
$$
\lim_{n\to\infty}\int_M |f_n(p) - g(p)|d\mu(p) = 0.
$$
In order to prove this, we prove that $f_n$ is a Cauchy sequence in $L^1(d\mu)$.
We take $\epsilon>0$ and choose $N$ large enough so that
\be \la{2.1}
\sup_{p\in M}\left |U_n(p) - \bar U(p)\right| \le \frac{\epsilon^2}{16 b},
\ee
and
$$
{\calE}_b[f_n] \le a + \frac{\epsilon^2}{16}
$$ 
holds for $n\ge N$. Let $s(p) = \frac{1}{2}(f_n(p) + f_m(p))$ with $n,m \ge N$. Then $\int_M sd\mu =1$,  $s\ge 0$, so
$$
a\le {\calE}_b[s].
$$
Therefore
$$
\frac{1}{2}\left \{{\calE}_b[f_n] + {\calE}_b[f_m]\right\} - {\calE}_b[s] \le \frac{\epsilon^2}{16}.
$$ 
On the other hand,
$$
\int_M \left \{\frac{1}{2}\left (f_n\log f_n + f_m\log f_m\right ) - s\log s\right\}d\mu \le 
\frac{1}{2}\left \{{\calE}_b[f_n] + {\calE}_b[f_m]\right\} - {\calE}_b[s] + \frac{\epsilon^2}{16}
$$
using \eqref{2.1} and $U[s]=\tfrac 12(U[f_n]+U[f_m])$, so
$$
\int_M \left \{\frac{1}{2}\left (f_n\log f_n + f_m\log f_m\right ) - s\log s\right\}d\mu \le \frac{\epsilon^2}{8}.
$$
Denote $\chi = \frac{f_n-f_m}{f_n + f_m}$ and note that $-1\le \chi\le 1$ holds 
$\mu$ - a.e. Also, elementary calculations show that
$$
\left \{\frac{1}{2}\left (f_n\log f_n + f_m\log f_m\right ) - s\log s\right\} =
\frac{s}{2} G(\chi)
$$
holds with
$$
G(\chi) = \log(1-\chi^2) + \chi\log\left(\frac{1+\chi}{1-\chi}\right).
$$
Note that $G$ is even on $(-1,1)$, that $G^{\prime}(\chi) = \log\left(\frac{ 1+\chi}{1-\chi}\right)$, $G(0) = G^{\prime}(0) = 0$ and $G^{\prime\prime}(\chi) = \frac{2}{1-\chi^2}\ge 2$ on $(-1,1)$. Consequently,
$$
0\le \chi^2 \le G(\chi)
$$
holds for $-1\le \chi \le 1$. It follows that we have 
$$
\int_M \frac{(f_n-f_m)^2}{f_n+f_m} d\mu \le \frac{\epsilon^2}{2}
$$
Writing $|f_n-f_m| = \sqrt{f_n+f_m}\frac{|f_n-f_m|}{\sqrt{f_n+f_m}}$
and using the Schwartz inequality we deduce
$$
\int_M|f_n-f_m|d\mu \le \epsilon.
$$
Therefore the sequence $f_n$ is Cauchy in $L^1(d\mu)$. This proves that
the weak limit $f_nd\mu \to gd\mu$ is actually strong $f_n\to g$ in $L^1(d\mu)$. By passing to a subsequence if necessary, we may assume that $f_n\to g$ holds also $\mu$- a.e. Then from Fatou's Lemma,
$$
\int_Mg\log gd\mu \le \lim_{j\to\infty}\int_M f_j\log f_j d\mu.
$$
and thus $g$ is a minimizer of $\calE_b$ with ${\mathcal{E}}_b[g] = a$.

It follows that $g\ge\delta$ where $\delta>0$ is such that $(x\log x)'<-3\|k\|_\infty$ for all $x\le\delta$. Otherwise $g_1=\max\{g,\delta\}$ satisfies $\calE_b[g_1]<\calE_b[g]$ and if we let $g_2=\min\{g_1,m\}$ with $m>1$ such that $\int_M g_2d\mu=1$, then $\calE_b[g_2]\le \calE_b[g_1] < \calE_b[g]$, a contradiction. Thus $g\ge\delta$ and then
the fact that $g$ solves the Onsager equation \eqref{geq} follows by taking the Gateaux derivative of $\calE_b$ with respect to each bounded $h\in L^1(M)$ with $\int_M h d\mu=0$. Since $U$ is  Lipschitz, the same is true about $g$. This concludes the proof.

\smallskip

The next two results describe the behavior of the minimizers of $\calE_b$ in the vanishing temperature limit $b\to\infty$. Let us assume from now on that there is no repulsion between identical (or identically oriented) particles, that is,
\be
k(p,p) = 0\quad \mbox{for all}\; p\in M.
\la{excl}
\ee

\begin{thm} \la{T.2}
Let $M$ be compact and metric, let $k$ satisfy (\ref{kprop}) and (\ref{excl}), and let $\mu$ be a Borel probability measure on $M$. Let $b_n\to\infty$ and consider a sequence of  free energy minimizers $g_n\ge 0$, $\int_M g_n d \mu =1$,
${\calE}_{b_n}[g_n] = \min_{\{f\ge 0, \; \int fd\mu =1\}}{\calE}_{b_n}[f]$. Then the sequence of measures $\nu_n = g_nd\mu$ has a subsequence that converges weakly to a Borel probability measure $\nu$ supported on a set $A\subseteq M$ such that $k(p,q)=0$ for any $p,q\in A$. 
\end{thm}

{\bf Remark.} Of course, this shows that all limit points of any sequence of minimizers of $\calE_{b_n}$ as $b_n\to\infty$ are as $\nu$ in the theorem. Notice that if $k(p,q)\neq 0$ for $p\neq q$, then each such $\nu$ must be $\delta_p$ for some $p\in M$.
\smallskip 

\noindent{\bf Proof.} There is obviously a subsequence of $\nu_n$ converging to a measure $\nu$. We relabel the subsequence $\nu_n$. The normalization $\nu_n(M) =1$ of the  probability measures $\nu_n$ shows that $\nu$ is a probability measure. 

We observe that for any $r>0$ there exists a ball $B = B(p,r)$ in $M$ such that $\mu(B)>0$. Indeed, if this is not the case, then there exists $r>0$ such that $\mu(B(p,r)) = 0$ for all $p\in M$.
Because $M$ is compact, we can cover it with finitely many such balls, and deduce $\mu(M) =0$, a contradiction.

Now we claim  that
\be \la{2.2}
\lim_{n\to\infty}\frac{2}{b_n}{\calE}_{b_n}[g_n] = 0.
\ee
Indeed, to prove this, we pick $\epsilon>0$ and use the uniform continuity of $k$ and the property (\ref{excl}) to find $r= r(\epsilon)$ so small that if $d(p,q)\le 2r$ then $k(p,q)\le \frac{\epsilon}{2}$.
We take a ball $B=B(p,r)$ such that $\mu(B)>0$ and consider the normalized indicator function of $B$, $f(p) = \mu(B)^{-1}{\mathbf 1}_{B}(p)$.   Then 
$$
\frac{2}{b_n}{\mathcal{E}}[g_n] \le \frac{2}{b_n}\log \mu(B)^{-1} + \frac{\epsilon}{2}
$$
holds because $g_n$ is an energy minimum. Fixing $r(\epsilon)$ we may find $N$ so that $\frac{2}{b_n}\log\left(\mu(B)\right)^{-1}\le \frac{\epsilon}{2}$  
holds for $n\ge N$ and thus \eqref{2.2} holds. 

Let now $A$ be the support of $\nu$ and assume there are $p,q\in A$ such that $k(p,q)=2\del>0$. Then for any $x\in B(p,\del/2L)$ and $y\in B(q,\del/2L)$ we have $k(x,y)\ge\del$ by \eqref{kprop}. Then 
$$
\liminf_{n\to\infty} \frac 2{b_n}\calE_{b_n}[g_n] \ge \del \nu(B(p,\tfrac \del{2L})) \nu(B(q,\tfrac \del{2L})) >0 ,
$$
contradicting \eqref{2.2}. The proof is finished.

\begin{thm} \la{T.3}
Let $M$ be compact and metric, let $k$ satisfy (\ref{kprop}) and (\ref{excl}), and let $\mu$ be a Borel probability measure on $M$.
Let $A_0,A_1\subseteq M$ be compacts with $k(p,q)=0$ for any $p,q\in A_j$ ($j=0,1$) and $B_j(\eps)=\{p\in M \,|\, d(p,A_j)<\eps\}$. Assume that for some $\eps_j>0$ there is a 1-1 map $T:B_1(\eps_1)\to B_0(\eps_0)$ (not necessarily onto) with $T$ and $T^{-1}$ measurable, and that there is $c>1$ such that 
\be \la{2.4}
\forall p,q\in B_1(\eps_1) \,:\, k(T(p),T(q))\le  k(p,q) 
\ee
\be \la{2.5}
\forall B\subseteq B_1(\eps_1) \text{ measurable }\,:\,  \mu(T(B))\ge c \mu (B). 
\ee
Assume also that for each $p\in A_1$, $q\in M\setminus B_1(\eps_1)$ we have $k(p,q)>0$.
Then $\nu (A_1)<1$ for each measure $\nu$ as in Theorem \ref{T.2}.
\end{thm}

{\bf Remark.} This and Theorem \ref{T.2} show that only those sets $A\subseteq M$ with $k(p,q)=0$ for $p,q\in A$ which have the largest (in the sense of $\mu$) neighborhoods can be supports of limit points of free energy minimizers. 
\smallskip

\noindent{\bf Proof.} 
Assume that the sequence of free energy minimizers $g_n\ge 0$ with $\int_M g_n d \mu =1$ and
$\calE_{b_n}[g_n] = \min_{\{f\ge 0, \; \int fd\mu =1\}}\calE_{b_n}[f]$ corresponds to 
$b=b_n$, the measures $\nu_n = g_nd\mu\rightharpoonup \nu$ weakly, and $\nu(A_1)=1$. We first claim $\mu(A_1)=0$. Indeed, otherwise, letting $A_2 = T(A_1)$ we have 
\[
\calE_{b_n}[\one_{A_2}] = -\log \mu(A_2) \le -\log \mu(A_1) -\log c. 
\]
But $\lim_{\eps\to 0} \mu(B_1(\eps))=\mu(A_1)$, $\lim_{n\to\infty} \int_{B_1(\eps)} g_n d\mu=1$ for each $\eps>0$, and  Jensen's inequality give 
\[
\liminf_{n\to\infty} \calE_{b_n}[g_n] \ge \liminf_{n\to\infty} \int_M g_n\log g_n d\mu \ge -\log \mu(A_1) > \calE_{b_n}[\one_{A_2}],
\]
a contradiction.

Without loss of generality assume that $L\ge 1$ in \eqref{kprop}.
Compactness and \eqref{kprop} imply that $\del_1= \inf\{k(p,q)\,|\, p\in A_1, \, q\in M\setminus B_1(\eps_1)\}>0$. Let $\gamma= \mu(B_1(\del/4L)\setminus B_1(\del/8L))\le 1$ where $\del\in (0,\min\{ \del_1,\eps_1\}]$ is chosen so that $\gamma>0$ (if such $\del$ does not exist, then $\mu(B_1(\min\{\del_1,\eps_1\}/4L))=0$ because $\mu(A_1)=0$, and so $\nu(A_1)=0$). 

For all large enough $n$ we have
\be \la{2.6}
\int_{B_1(\tfrac \del{8L})} g_n d\mu \ge 1-\omega
\ee
where $\omega>0$ will be specified later. We also let 
\be \la{2.7}
\alpha_n= \int_{M\setminus B_1(\eps_1)} g_n d\mu \le\omega
\ee
 (because $\del/8L\le\del\le \eps_1$).

Let $\mu_* = T_*(\mu|_{B_1(\eps_1)})$ be the pushforward measure of $\mu$ restricted to $B_1(\eps_1)$ (thus $\mu_*$ is supported in $T(B_1(\eps_1))\subseteq B_0(\eps_0)$). Then \eqref{2.5} shows that $\mu_*$ is absolutely continuous with respect to $\mu$ and that there is a measurable function $0\le f\le c^{-1} \one_{T(B_1(\eps_1))}$ such that $d\mu_*=f d\mu$. Define 
\[
 f_n(p) = \begin{cases} g_n(T^{-1}(p))f(p) & p\in T(B_1(\eps_1)), \\ 0 & p\in M\setminus T(B_1(\eps_1)) \end{cases}
\]
and notice that for any $B\subseteq T(B_1(\eps_1))$,
\be \la{2.8}
\int_B f_n d\mu=\int_{T^{-1}(B)} g_n d\mu,
\ee
in particular, $\int_M f_n d\mu = 1-\alpha_n$. We let 
\[
\Sigma_n= \{p\in T(B_1(\tfrac \del{4L})\setminus B_1(\tfrac\del{8L})) \,|\, f_n(p)\le \tfrac 12\}
\]
so that $\mu(\Sigma_n)\ge \tfrac \gamma 2$ by \eqref{2.8}, \eqref{2.6}, \eqref{2.5}, provided we take  $\omega\le \tfrac \gamma 4$. If we now let $\beta_n= \alpha_n/\mu(\Sigma_n)$, then $\beta_n\le \tfrac 12$ by \eqref{2.7} and $\omega\le \tfrac \gamma 4$, and $h_n= f_n + \beta_n \one_{\Sigma_n}$ satisfies $\int_M h_n d\mu =1$. We will now show that if $\omega$ is small enough and $n$ large, then $\calE_{b_n}[h_n]<\calE_{b_n}[g_n]$, thereby obtaining a contradiction.

For large $n$ we have by \eqref{2.4}, \eqref{2.6}, \eqref{2.7}, \eqref{2.8}, and \eqref{kprop},
\[
 \int_{B_1(\eps_1)^2} k(p,q) g_n(p)g_n(q) d\mu(p)d\mu(q) \\
\ge \! \int_{T(B_1(\eps_1))^2} k(p,q) f_n(p)f_n(q) d\mu(p)d\mu(q),
\]
\begin{align*}
\int_{B_1(\tfrac \del{8L})\times(M\setminus B_1(\eps_1))} k(p,q) g_n(p)g_n(q) d\mu(p)d\mu(q) & \ge (\del_1- L\tfrac \del {8L}) (1-\omega)\alpha_n \\
& \ge \tfrac {7\del}8  (1-\omega)\alpha_n,
\end{align*}
\[
\int_{T(B_1( \tfrac \del{8L}))\times \Sigma_n} k(p,q) f_n(p)\beta_n d\mu(p)d\mu(q) \\
\le L(\tfrac \del{8L}+\tfrac \del{4L}) \alpha_n \le \tfrac {3\del}8 \alpha_n,
\]
\[
 \int_{(M\setminus T(B_1(\tfrac \del{8L})))\times \Sigma_n} k(p,q) f_n(p)\beta_n d\mu(p)d\mu(q) \le \|k\|_\infty \omega \alpha_n,
\]
\[
 \int_{ \Sigma_n^2 } k(p,q) \beta_n\beta_n d\mu(p)d\mu(q) \le \|k\|_\infty  \alpha_n^2 \le \|k\|_\infty \omega \alpha_n,
\]
where in the third line we have used
\[
 k(T(p),T(q))\le k(p,q)\le k(p',q')+L(d(p,p')+d(q,q')) = L(d(p,p')+d(q,q'))
\]
for $p',q'\in A_1$. Thus \eqref{2.7} gives
\be \la{2.9}
\int_{M^2} k(p,q) g_n(p)g_n(q) d\mu(p)d\mu(q) > \int_{M^2} k(p,q) h_n(p)h_n(q) d\mu(p)d\mu(q)
\ee
provided $\tfrac{7\del}4 (1-\omega)> \tfrac{3\del}4 + 3\|k\|_\infty \omega$ (or $\omega \le \del(8\|k\|_\infty)^{-1}$ since $\omega\le \tfrac\gamma 4\le \tfrac 14$).
The definition of $\Sigma_n$, $\beta_n\le \tfrac 12$, $(x\log x)'\le 1$ for $x\le 1$, $f\le c^{-1}$, \eqref{2.8}, and \eqref{2.7} give
\begin{align*}
\int_M h_n\log h_n d\mu &\le \int_M f_n\log f_n d\mu + \alpha_n \le \int_{B_1(\eps_1)} g_n \log(g_nc^{-1}) d\mu +\alpha_n \\
&\le  \int_{B_1(\eps_1)} g_n \log g_n d\mu - (1-\omega)\log c +\alpha_n,
\end{align*}
while Jensen's inequality and \eqref{2.7} give
\[
 \int_{M\setminus B_1(\eps_1)} g_n \log g_n d\mu \ge \alpha_n \log \frac{\alpha_n}{\mu(M\setminus B_1(\eps_1))} \ge \alpha_n \log \alpha_n.
\]
If now $ \alpha_n (1-\log\alpha_n) < (1-\omega)\log c$ (which is guaranteed if $\omega$ is such that $\omega(1-\log\omega)< \tfrac 34 \log c$ because $\alpha_n\le\omega\le \tfrac 14$), then
\[
 \int_M h_n\log h_n d\mu < \int_M g_n\log g_n d\mu
\]
and we have $\calE_{b_n}[h_n]<\calE_{b_n}[g_n]$ for large $n$. Thus we only need to pick 
\[
 \omega< \min \left\{ \frac \gamma 4, \frac \del{8\|k\|_\infty} \right\} \quad\text{and}\quad \omega(1-\log\omega) < \frac {3}4 \log c 
\]
and the proof is finished.

\section{Examples}

Examples of rod-like particles have been discussed in detail \cite{c1,ckt} and will not be discussed here. It suffices to say that Theorem {\ref{T.2}}
implies that the zero temperature limit of the minimizers has to be a prolate state, both for the Maier-Saupe and the Onsager potentials, because $k(p,q) \neq 0$ for $p\neq q$ in both cases.
\smallskip

{\bf Example 1.} We start the list of examples with interacting two-rods, each made of two segments of unit length attached at the origin. The corpora belong to $M = {\mathbb S}^1\times {\mathbb S}^1$. We assume that the interaction between pairs of two-rods is determined entirely by the {\it area and orientation} of the triangle formed by each two-rod, and each two-rod corpus rejects two-rods that have very different oriented area than itself. The simplest interaction that achieves this is
$$
k((p_1,p_2), (q_1,q_2)) =  \|e(p_1)\wedge e(p_2) -e(q_1)\wedge e(q_2)\|^2
$$
with $e(p) = (\cos p, \sin p)$ if $p \in [0, 2\pi)$. Each ${\mathbb{S}}^1$
is viewed as a subset of ${\mathbb R}^2$; the exterior product $\bigwedge^2 ({\mathbb R}^2)$ is isomorphic to ${\mathbb R}$. Every element in it is a multiple of $e_1\wedge e_2$, with $e_1 = (1,0)$ and $e_2 = (0,1)$, so $e(p_1)\wedge e(p_2) = \sin(p_1-p_2)(e_1\wedge e_2)$. We then have
$$
k((p_1,p_2), (q_1,q_2)) = \left (\sin(p_1-p_2) - \sin(q_1-q_2)\right)^2.
$$
We take the uniform measure $d\mu(p_1,p_2) = \frac{1}{4\pi^2}dp_1 dp_2$ on ${\mathbb S}^1\times{\mathbb S}^1$ and note that the potential has the form
$$
U[f](p_1, p_2) =  \sin^2(p_1-p_2) -2a\sin(p_1-p_2) + \gamma
$$
with 
\be \la{3.0}
 a=\int_M \sin(q_1-q_2) f(q_1,q_2) d\mu,
\ee
\[
\gamma = \int_M \sin^2(q_1-q_2) f(q_1,q_2) d\mu.
\]
Onsager's equation $f = (Z_b[f])^{-1}e^{-bU[f]}$ now gives 
\[
f(p_1,p_2)=Z^{-1}e^{-b(\sin(p_1-p_2)-a)^2}
\]
with 
\[
Z=Z_b[f]e^{b(\gamma-a^2)} = \frac 1{2\pi} \int_0^{2\pi}e^{-b(\sin \theta  - a)^2}d\theta.
\]
Thus by \eqref{3.0}, solving it is equivalent to finding $a$ such that
\be \la{3.0a}
a = [\sin\theta]_{b} (a)
\ee
where
$$
[\phi]_{b} (a) = \frac{\int_0^{2\pi}\phi (\theta)  e^{-b(\sin \theta - a)^2}d\theta}{\int_0^{2\pi} e^{-b(\sin \theta  - a)^2}d\theta}.
$$

We note that $a=0$ is always a solution that yields
$$
f_b(p_1,p_2) = Z^{-1}e^{-b\sin^2(p_1-p_2)}.
$$
As $b\to\infty$ this tends to $\delta ((p_1-p_2)\, \mbox{mod}\,\pi)$, the uniform measure on the union of the segments $p_1-p_2=0$ and $p_2-p_2=\pi$, whose support corresponds to two-rods with area 0. However, this is not the only solution for large $b$.


Consider
$$
h_b(a) = \int_0^{2\pi} u e^{-bu^2}d\theta
$$
with $u (\theta, a) = \sin\theta -a$. In order to solve Onsager's equation (i.e., (\ref{3.0a})) we seek  zeros of $h_b(a)$. Clearly 
$h_b(1)<0$. In order to prove that $h_b$ vanishes for some positive $a$, it is therefore enough to find $0<a<1$ such that $h_b(a)>0$. Now, changing variables, we have
$$
h_b(a) = 2\int_{-1-a}^{1-a} e^{-bu^2}\frac{udu}{\sqrt{1-(u+a)^2}}
$$
and fixing $0<a<1$ we get
\[
\ba
h_b(a)  =  4\int_0^{1-a}e^{-bu^2}u\left\{\frac{1}{\sqrt{1-(u+a)^2}} - \frac{1}{\sqrt{1-(a-u)^2}}\right\} \!du 
 -2\int_{1-a}^{1+a}e^{-bu^2}\frac{udu}{\sqrt{1-(a-u)^2}} \\
 =  16a\int_0^{1-a}e^{-bu^2}u^2\frac{du}{\sqrt{(1-(u+a)^2)(1-(u-a)^2)}\left(\sqrt{1-(u-a)^2}+ \sqrt{1-(u+a)^2}\right)}\\ 
 -2\int_{1-a}^{1+a}e^{-bu^2}\frac{udu}{\sqrt{1-(a-u)^2}}.
\ea
\]
The last integral is negative but exponentially small as $b\to\infty$.
The positive integral is larger than a fixed multiple of $\int_0^\delta u^2e^{-bu^2}du$ for small fixed $\delta$ (depending on $a$). This integral asymptotically equals $Cb^{-3/2}$ with positive $C$, as $b\to\infty$. Thus, for any fixed $0<a<1$ and $b$ large enough, we have $h_b(a)>0$. This proves for all large enough $b$ the existence of $a(b)>0$ such that $h_b(a(b))= 0$. Moreover, because $0<a<1$ was arbitrary, this also proves  $\lim_{b\to\infty} a(b) = 1$ for any such $a(b)>0$. A similar argument applies for $0>a>-1$, and consequently limit points of solutions $f_b$ of Onsager's equation as $b\to\infty$ are precisely $\delta ((p_1-p_2)\, \mbox{mod}\,\pi)$, $\delta \left((p_1-p_2 -\frac{\pi}{2} )\, \mbox{mod}\,2\pi \right)$, and $\delta\left ((p_1-p_2 + \frac{\pi}{2} )\, \mbox{mod}\,2\pi \right)$.

We will now employ Theorem \ref{T.3} to show that the first of these cannot be the limit point of minimizers of the free energy as $b\to \infty$. Assume it is, let $A_1=\{(p_1,p_2)\,|\, \sin(p_1-p_2)=0 \}$ and $A_0=\{(p_1,p_2)\,|\, \sin(p_1-p_2)=1 \}$. For some $c>1$ consider the map 
\[
T(p_1,p_2) = \begin{cases} (p_2+\tfrac \pi 2 -c(\eps-p_1+p_2) ,p_2) & p_1-p_2\in(-\eps,\eps) \\
              (p_2+\tfrac \pi 2 - c( \pi-\eps -p_1+p_2) ,p_2) & p_1-p_2\in (\pi-\eps,\pi+\eps)
\end{cases} 
\]
from $B_1(\eps)=\{(p_1,p_2)\,|\, p_1-p_2\in(-\eps,\eps)\cup(\pi-\eps,\pi+\eps) \}$ onto $B_0(2c\eps)=\{(p_1,p_2)\,|\, p_1-p_2\in(\tfrac \pi 2-2c\eps,\tfrac \pi 2+2c\eps) \}$. If $\eps>0$ is small enough, then the conditions of Theorem \ref{T.3} are satisfied, with \eqref{2.4} holding because 
\[
 (\sin x- \sin(x+\tilde\eps))^2 = \tilde\eps^2 \cos^2 x + O(\tilde\eps^3) 
\]
and $\cos^2$ is continuous with $\cos^2 \tfrac \pi 2=0$ and $\cos^2 0=\cos^2\pi=1$.
The conclusion of the theorem then contradicts our hypothesis.

This and symmetry show that the limit points of minimizers of the free energy as $b\to \infty$ are precisely the two measures $\delta\left ((p_1-p_2 \pm \frac{\pi}{2} )\, \mbox{mod}\,2\pi \right)$, each supported on the set of two-rods forming triangles with area $\tfrac 12$ but with opposite orientations. 

\smallskip

{\bf Remark.} We note that Theorem \ref{T.2} a priori shows that any  limit point of free energy minimizers  as $b\to\infty$ must be a measure supported on some set $M_a=\{(p_1,p_2)\,|\, \sin(p_1-p_2)=a \}$. A variant of the above argument involving Theorem \ref{T.3} then excludes any $a\neq\pm 1$. However, it does not show that the limiting measures are uniform on $M_{\pm 1}$ (that follows from symmetry), nor does it capture the existence of the solutions $f_b$ of Onsager's equation.

\smallskip

{\bf Example 2.} Let us now consider an example in which we allow the sizes of the 
rods to vary. We parameterize the corpora by $(x_1, x_2, p_1,p_2)\in M=[0,L]^2\times  [0,2\pi]^2$ and
let $d\mu = \frac{1}{4\pi^2 L^2}dxdydp_1dp_2$.  Each corpus is a two-rod with
segments of lengths $x_1, x_2\in[0,L]$ emanating from the origin at angles $p_1,p_2\in [0,2\pi]$. We consider
\begin{align*}
& U [f] (x_1,x_2,p_1,p_2) \\
& = \int_M (x_1x_2\sin(p_1-p_2)-y_1y_2\sin(q_1-q_2))^2f(y_1,y_2,p_1,p_2)d\mu(y_1,y_2,q_1,q_2). 
\end{align*}
It is again easy to see that solutions of Onsager's equation are of the form
$$
f(x_1,x_2,p_1,p_2) = \frac {e^{-b(x_1x_2\sin(p_1-p_2) -a)^2}} {\int_M e^{-b(x_1x_2\sin(p_1-p_2) -a)^2} d\mu}
$$
with $a$ determined by
$$
a= \int_M x_1x_2\sin(p_1-p_2)\, f(x_1,x_2,p_1,p_2)d\mu
$$
Introducing, as above, the function 
$$
u(x_1,x_2,\theta,  a) = x_1x_2\sin\theta -a
$$
and the notation
$$
[\phi]_b(a) =\frac {\int_{M'} \phi e^{-b(x_1x_2\sin\theta -a)^2} d\mu'}  {\int_{M'} e^{-b(x_1x_2\sin\theta -a)^2} d\mu'}
$$
with $M'=[0,L]^2\times[0,2\pi]$ and $d\mu' = \frac{1}{2\pi L^2}dxdyd\theta$, we see that $a$ is determined by the requirement $[u]_b(a)= 0$. 


We have $[u]_b(a)=0$ if and only if 
$$
\int_{M'} ue^{-bu^2}d\mu' = 0.
$$
Now, in view of
$$
\frac{\partial u}{\partial x_2} = x_1\sin\theta
$$
we can write
$$
\int_{M'} ue^{-bu^2}d\mu' = \frac{1}{4\pi b L^2}\int_0^{2\pi}\int_0^L\frac{1}{x\sin\theta}\left(e^{-ba^2}-e^{-b(Lx\sin\theta-a)^2}\right)dxd\theta
$$
Thus $a$ must obey
$$
0 = \int_0^{2\pi}\int_0^L\left(e^{bLx\sin\theta(2a-Lx\sin\theta)} - 1\right)\frac{dx d\theta}{x\sin\theta},
$$
where the integrand is bounded. 

We will now show that if $a_n$ solve $[u]_{b_n}(a_n)=0$ for some sequence $b_n\to\infty$, then $a_n\to 0$. 
Assume $a>0$, so then clearly
$$
I_1= \int_{Lx\sin\theta\ge 2a}\left(e^{bLx\sin\theta(2a-Lx\sin\theta)} - 1\right)\frac{dx d\theta}{x\sin\theta}
$$
obeys $-\frac{\pi L}{2a}\le I_1 \le 0$ and the integrand in 
$$
I_2 = \int_{Lx\sin\theta\le 2a}\left(e^{bLx\sin\theta(2a-Lx\sin\theta)} - 1\right)\frac{dx d\theta}{x\sin\theta} 
$$
is positive. Assume that some subsequence of $a_n$, which we again call $a_n$, obeys $a_n\ge \alpha>0$ for some fixed $\alpha>0$. Then
$$
I_2 \ge \! \! \int_{Lx\sin\theta\le \alpha}\left(e^{b_n\alpha Lx\sin\theta}-1\right)\frac{dxd\theta}{x\sin\theta}\ge b_n\alpha L \left| \left\{(x,\theta) \,|\, 0\le x\sin\theta \le \alpha\right\} \right|\to \! \infty
$$
while $I_1$ remains bounded. Thus $I_1+I_2 >0$ for large $n$, a contradiction. A similar argument excludes a subsequence with $a_n\le -\alpha<0$, proving $a_n\to 0$.

Thus, in contrast to Example 1, this time all solutions of Onsager's equation (and therefore also minimizers of the free energy) converge to a measure supported on the set of two-rods forming triangles with zero area as $b\to \infty$.

\smallskip

{\bf Example 3.} Here we demonstrate another application of Theorem \ref{T.3}. This time our corpora are rhombi  $R_p$  of unit side, lying in $\mathbb R^2$, aligned with each other --- centered at the origin, their longer diagonal coinciding with the $x$-axis. They are parameterized by the smaller angle $p\in[0,\tfrac \pi 2]=M$ and $d\mu = \tfrac 2\pi dp$ is the uniform measure. The kernel $k$ equals the area of the {\it symmetric difference} of the rhombi so that a rhombus rejects those with which it has small overlap. We then get
\be \la{3.1}
 k(p,q)= 
8\frac {\sin^2 \frac {|p-q|}4 } {\sin \frac {|p-q|}2 } 
\left( \sin^2 \frac {p+q}4 \sin \frac {p}2 \sin \frac {q}2  + \cos^2 \frac {p+q}4 \cos \frac {p}2 \cos \frac {q}2 \right),
\ee
which is obtained as follows. Let the side of $R_p$ lying in the first quadrant be bisected by that of $R_q$ into segments of lengths $a,1-a$ (the former having one end at the $x$-axis). Similarly, the side of $R_q$ is bisected into segments of lengths $b,1-b$. Then we obtain $a \sin \tfrac p 2 = b\sin \tfrac q 2$ and $(1-a)\cos \tfrac p 2 = (1-b) \cos \tfrac q 2$. Evaluating $a, b$ from this and then computing the area of the symmetric difference gives \eqref{3.1}.

Assume that $g_n$ is a sequence of free energy minimizers for some $b_n\to\infty$. Theorem \ref{T.2} shows that a subsequence of the measures $g_nd\mu$ converges to the measure $\nu=\delta_q$ for some $q\in[0,\tfrac \pi 2]$ (each set $A$ satisfying the condition in Theorem \ref{T.2} contains a single element so, any limit $\nu$ must be some $\delta_q$). The next question is which $\delta_q$ can be limiting points of minimizers of the free energy.

Note that
\be \la{3.2}
 k(p,p+\eps) = \eps \left( \sin^4 \frac p 2 + \cos^4 \frac p 2 \right) + O(\eps^2),
\ee
with the bracket being smallest (equal to $\tfrac 12$) when $p=\tfrac \pi 2$. If $q\neq \tfrac \pi 2$, let $A_0=\{\tfrac \pi 2\}$, $A_1= \{q\}$, and consider the map $T(p)=\tfrac \pi 2 - c(q+\eps-p)$ from $B_1(\eps)= (q-\eps,q+\eps)$ to $B_0(2c\eps)= (\tfrac \pi 2 - 2c\eps, \tfrac \pi 2)$. Here $c\in \left(1 , \sin^4 \frac q 2 + \cos^4 \frac q 2 + \tfrac 12 \right)$ and $0<\eps\ll \tfrac \pi 2 - q$. If $\eps$ is small enough, we find using \eqref{3.2} that the assumptions of Theorem \ref{T.3} hold and thus $\nu(\{q\})<1$, a contradiction. Therefore $q=\tfrac \pi 2$ for any $b_n\to\infty$ and any convergent subsequence of $g_n d\mu$ so, in fact, if $g_b$ is a minimizer of $\calE_b$, then we have the weak convergence $g_b d\mu \rightharpoonup \delta_{\pi/2}$ to the delta function on the square $R_{\pi/2}$ as $b\to\infty$.

\section{Conclusions and Outlook}
The zero temperature limit of interacting corpora, under the influence of disorder and conformation constraints is supported by ur-corpora. These are selected by an entropic popularity contest: the winners are corpora $p\in M$ with most (in the sense of measure) conforming corpora $q$ (with $k(p,q) \approx \min k = 0$).
 
Possible extensions of the theory to non-compact $M$ are of interest. A kinetic theory \cite{doied} exists in the case of a compact Riemannian manifold and $\mu$ the normalized volume measure \cite{ckt,ctv}. The natural extension of this theory, in the spirit of \cite{ags} is being pursued. The coupling of this to macroscopic fluids, in the spirit of \cite{c-nfp,c-smo,cftz,cm}, is a further goal.
\smallskip

{\bf{Acknowledgment}} The research of PC was partially supported by the NSF grant
DMS-0504213.  The research of AZ was partially supported by the NSF grant DMS-0632442 and an Alfred P. Sloan Research Fellowship.

\end{document}